\documentclass[preprint,10pt]{elsarticle}
\usepackage{amsmath}
\usepackage{amsthm}
\usepackage{amssymb}
\usepackage{algorithm}
\usepackage{algpseudocode} 
\usepackage{url}
\usepackage{enumerate}
\usepackage[margin=1.25in]{geometry}
\usepackage{xcolor}

\newtheorem{theorem}{Theorem}[section]
\newtheorem{lemma}[theorem]{Lemma}
\newtheorem{proposition}[theorem]{Proposition}
\newtheorem{corollary}[theorem]{Corollary}

\newtheorem{remark}[theorem]{Remark}

\newcommand{\qbinom}[2]{\genfrac{[}{]}{0pt}{}{#1}{#2}}

\begin{document}

\begin{frontmatter}
\title{An evaluation algorithm for $q$--B\'ezier
triangular patches formed by convex combinations \tnoteref{fund}}

\tnotetext[fund]{ This work was partially supported through the Spanish research grant PGC2018-096321-B-I00 (MCIU/AEI) and by Gobierno de Arag\'on (E41\_20R).}

\author[rvt]{J. Delgado\corref{cor1}}
\ead{jorgedel@unizar.es}

\author[rvt]{H. Orera}
\ead{hectororera@unizar.es}

\author[rvt]{J.~M. Pe\~{n}a}
\ead{jmpena@unizar.es}

\cortext[cor1]{Corresponding author}

\address[rvt]{Departamento de Matem\'{a}tica Aplicada/IUMA, Universidad de Zaragoza, Spain}

\begin{abstract}
An extension to  triangular domains of the univariate $q$-Bernstein basis functions is introduced and analyzed. Some recurrence relations and properties such as partition of unity and degree elevation are proved for them. It is also proved that they form a basis for the space of polynomials of total degree less than or equal to $n$ on a triangle. In addition, it is presented a de Casteljau type evaluation algorithm whose steps are all linear convex combinations.
 
\end{abstract}
\end{frontmatter}

\section{Introduction}

Efficient evaluation of curves and surfaces is an important issue in computational mathematics 
(see \cite{MP,SIAM,ACM}). 
The $q$--Bernstein ($0<q\le 1$) basis of univariate polynomials has played an important role in several fields, such as Computer Aided Geometric Design (CAGD), Approximation Theory or Quantum Calculus. They have received much attention in recent research  (cf. \cite{gs}, \cite{sd}, \cite{LW} and
references in there). For the particular case $q=1$, it coincides with the basis of Bernstein polynomials. 

Bernstein polynomials were introduced more than a century ago and they later had a key role in CAGD \cite{f}. The Bernstein-B\'ezier model for the design of curves and surfaces provides basic algorithms (see \cite{farin}) and shape preserving properties (see \cite{cp1}) in this area. Besides, they are also very useful in other areas such as Numerical Analysis, Approximation Theory or, more recently, Partial Differential Equations. By means of barycentric coordinates, B\'ezier curves can be generalized to surfaces, yielding triangular patches, which provide a more natural generalization of B\'ezier curves than that provided by tensor product patches.

Corner cutting algorithms form a fundamental family of algorithms in CAGD. The basic steps of these algorithms are formed by linear convex combinations. In addition to their geometric interpretation, they present nice stability properties \cite{MP}. In fact, the well-known de Casteljau algorithm for the evaluation of polynomial curves is an example of corner cutting algorithm.

Among other generalizations of B\'ezier curves, $q$--B\'ezier curves provide more versatility to classical B\'ezier curves due to the use of the parameter $q$. Algorithms and geometric properties have been studied (see \cite{gs}, \cite{sd}, \cite{LWAG}, \cite{qBer}, \cite{cast}, \cite{Ph1}, \cite{qBlos} and references in there). In \cite{dp1}, algorithms for solving algebraic problems for collocation matrices of $q$--Bernstein polynomials with high relative accuracy were presented. Rational $q$--Bernstein bases with their properties were analyzed in
\cite{RQBBC}. Tensor product and rational $q$--B\'ezier  surfaces were studied in \cite{TPQB} and \cite{dp2}. For another generalized Bernstein basis, the $h$-Bernstein basis (see \cite{gs,hBer}), the extension to a basis over a triangular domain has been recently obtained (see \cite{LLRS}). The $q$--Bernstein  basis of univariate polynomials has been generalized to a triangular domain in \cite{LWAG}. However, the evaluation algorithm proposed in \cite{LWAG} is not in general formed by linear convex combinations. So, the usefulness of  evaluation algorithms based on convex combinations motivated us to research
for new extensions that always guarantee their existence. In this paper, we propose an alternative approach of $q$-Bernstein polynomials over a triangular domain that always provides such an evaluation algorithm. Moreover, the particular univariate application of our approach also has a corner cutting evaluation algorithm, as well as the natural corresponding extension to tensor products of $q$-Bernstein bases.

The paper is organized as follows. In Section 2, we present basic notations. In Section 3, we present the triangular $q$--Bernstein polynomials and we prove that they satisfy two recurrence relations. We also provide graphical examples of 
triangular $q$--Bernstein  polynomials, illustrating the effect of changing the parameter $q$. In Section 4, we obtain a de Casteljau type evaluation algorithm formed by linear convex combinations for $q$--B\'ezier representations of bivariate polynomials and the corresponding univariate evaluation algorithm. We also prove that the triangular $q$--Bernstein  polynomials form a partition of the unity. In Section 5, we present a degree elevation formula for triangular $q$--Bernstein  polynomials and we prove that the triangular $q$--Bernstein  polynomials of degree $n$ form a basis of the space of bivariate polynomials of degree less than or equal to $n$. We define $q$--B\'ezier patches and we include graphical examples of $q$--B\'ezier patches, illustrating again the effect of changing the parameter $q$. Although the basis of $q$-Bernstein polynomials over triangles has a stable evaluation algorithm because all its steps are linear convex combinations, in Section 6 we shall see that the usual Bernstein basis over triangles is always better conditioned. Finally, Section 7 summarizes the main conclusions of the paper.

\section{Basic notations}

Let us consider a nondegenerate triangle $\mathcal{T}:=\langle T_1,T_2,T_3\rangle$ on $\mathbb{R}^2$ whose vertices are the three points $T_i:=(x_i,y_i)$ for $i=1,2,3$. Then we know that every point $P:=(x,y)\in\mathbb{R}^2$ has a unique representation in the form 
\begin{equation}\label{coord}
P=uT_1 +v T_2 + w T_3,
\end{equation}
with $1=u+v+w$. The numbers $u, v$ and $w=1-u-v$ are called the \textit{barycentric coordinates} of the point $P$ relative to the triangle $\mathcal{T}$. If the three coordinates are positive, then the point lies in the interior of $\mathcal{T}$.

Given a positive real number $q$ and a natural number $r$ we define the {\it $q$-integer} $[r]$ as 

\[[r]:=\left\{\begin{array}{ll}
1 +q +\cdots + q^{r-1}=\frac{1-q^r}{1-q}, & \text{ if } q \neq 1, \\ 
r, & \text{ if } q=1.
\end{array}\right.\]
Then we can define in terms of the $q$-integers the following $q$-analogues. 
The {\it $q$-factorial} $[r]!$ (see \cite{kc}) is given by
\[[r]!:=\left\{\begin{array}{ll}
[r][r-1]\cdots [1], & \text{ if } q \neq 1, \\ 
r!, & \text{ if } q=1,
\end{array}\right.\]
and the {\it $q$-binomial coefficient}  $\qbinom{i}{j}$ is defined as
\begin{equation}\label{def.qbi}
\qbinom{i}{j}:=\frac{[i]!}{[j]![i-j]!}
\end{equation}
if $i\geq j \geq 0$ and as $0$ otherwise.
Let us recall the recurrence relation that defines the classical binomial coefficients $\binom{n}{k}$,
\begin{equation}\label{rec.binom}
\binom{n}{k}=\binom{n-1}{k}+\binom{n-1}{k-1}.
\end{equation}
The $q$-binomial coefficients satisfy the following recurrence relations, which are $q$-analogues of \eqref{rec.binom}:
\begin{equation}\label{rec1}
\qbinom{i}{j}=\qbinom{i-1}{j-1} + q^j \qbinom{i-1}{j},
\end{equation}
\begin{equation}\label{rec2}
\qbinom{i}{j}=q^{i-j}\qbinom{i-1}{j-1} +  \qbinom{i-1}{j}.
\end{equation}

\section{An extension of the triangular Bernstein polynomials}

Given a nonnegative integer $n$ and three indices $i,j,k\geq 0$ such that $i+j+k=n$, we define the following $q$-analogues of the triangular Bernstein polynomials relative to $\mathcal{T}$:

\begin{equation}\label{def.qBer}
B^n_{ijk}(u,v)=\qbinom{n}{k}\binom{i+j}{i}u^iv^j\prod_{s=0}^{k-1}(1-q^su-q^sv), \quad i,j,k\geq 0, \quad i+j+k=n.
\end{equation}
 We will call these polynomials the  \textit{triangular} $q$--\textit{Bernstein polynomials}. Let us notice that definition \eqref{def.qBer} coincides with the triangular Bernstein polynomials for $q=1$, and that the functions are nonnegative for all $q\in (0,1]$. So, from now on, we assume that $q\in(0,1]$. The first important property about the triangular $q$--Bernstein polynomials is that they can be defined in terms of the following recurrence relationships. 

\begin{proposition}
The triangular $q$--Bernstein polynomials \eqref{def.qBer} are given by 
\begin{equation}\label{rec}
B_{ijk}^n(u,v)=u B_{i-1,j,k}^{n-1}(u,v) + vB_{i,j-1,k}^{n-1}(u,v) + (q^{i+j}-q^{n-1}u-q^{n-1}v)B_{i,j,k-1}^{n-1}(u,v)
\end{equation}
and
\begin{equation}\label{recq}
B_{ijk}^n(u,v)=q^k u B_{i-1,j,k}^{n-1}(u,v) + q^k vB_{i,j-1,k}^{n-1}(u,v) + (1-q^{k-1}u-q^{k-1}v)B_{i,j,k-1}^{n-1}(u,v),
\end{equation}
where  $B_{000}^0(u,v)=1$ and $B_{ijk}^n(u,v)=0$ for any negative index $i,j,k$.
\end{proposition}
\begin{proof}
 Given a nonnegative integer $n$ and the nonnegative indices $i+j+k=n$, let us compute the polynomial  $p_{ijk}(u,v)$ defined by \eqref{rec}
\begin{align*}
p_{ijk}(u,v)=&u B_{i-1,j,k}^{n-1}(u,v) + vB_{i,j-1,k}^{n-1}(u,v) + (q^{i+j}-q^{n-1}u-q^{n-1}v)B_{i,j,k-1}^{n-1}(u,v)\\
=&\qbinom{n-1}{k} u^iv^j\prod_{s=0}^{k-1}(1-q^su-q^sv)\left( \binom{i+j-1}{i-1} +  \binom{i+j-1}{i} \right)\\
&+ \qbinom{n-1}{k-1}\binom{i+j}{i}  u^iv^j\prod_{s=0}^{k-2}(1-q^su-q^sv)(q^{i+j}-q^{n-1}u-q^{n-1}v).
\end{align*}
Taking into account that $q^{i+j} -q^{n-1}u-q^{n-1}v=q^{n-k}(1-q^{k-1}u-q^{k-1}v)$ and using \eqref{rec.binom}, we can rewrite $p_{ijk}(u,v)$ as
\begin{align*}
p_{ijk}(u,v)=& \binom{i+j}{i} u^iv^j\prod_{s=0}^{k-1}(1-q^su-q^sv)\left(\qbinom{n-1}{k} + q^{n-k}\qbinom{n-1}{k-1}\right),
\end{align*}
and by \eqref{rec2} we see that 
\begin{align*}
p_{ijk}(u,v)=& \qbinom{n}{k}\binom{i+j}{i} u^iv^j\prod_{s=0}^{k-1}(1-q^su-q^sv).
\end{align*}
Hence, $p_{ijk}(u,v)=B_{ijk}^n(u,v)$ and we conclude that \eqref{rec} gives a recursive definition of the triangular $q$--Bernstein polynomials \eqref{def.qBer}.

Let us now compute the polynomial $r_{ijk}(u,v)$ defined by the recurrence \eqref{recq},
\begin{align*}
r_{ijk}(u,v)=&q^k u B_{i-1,j,k}^{n-1}(u,v) + q^k vB_{i,j-1,k}^{n-1}(u,v) + (1-q^{k-1}u-q^{k-1}v)B_{i,j,k-1}^{n-1}(u,v)\\
=&q^k\qbinom{n-1}{k} u^iv^j\prod_{s=0}^{k-1}(1-q^su-q^sv)\left( \binom{i+j-1}{i-1} +  \binom{i+j-1}{i} \right)\\
&+ (1-q^{k-1}u-q^{k-1}v)\qbinom{n-1}{k-1}\binom{i+j}{i}  u^iv^j\prod_{s=0}^{k-2}(1-q^su-q^sv).
\end{align*}
Using \eqref{rec.binom} we deduce that  
\begin{align*}
r_{ijk}(u,v)=& \binom{i+j}{i} u^iv^j\prod_{s=0}^{k-1}(1-q^su-q^sv)\left(q^k\qbinom{n-1}{k} + \qbinom{n-1}{k-1}\right),
\end{align*}
and by \eqref{rec1} we have that
\[r_{ijk}(u,v)=\qbinom{n}{k}\binom{i+j}{i} u^iv^j\prod_{s=0}^{k-1}(1-q^su-q^sv).\]
And so, we conclude that $r_{ijk}(u,v)=B_{ijk}^n(u,v)$ and that \eqref{recq} gives a second recursive definition for the triangular $q$--Bernstein polynomials.
\end{proof}

\begin{figure}[!h]
\begin{center}
\includegraphics[width=\textwidth]{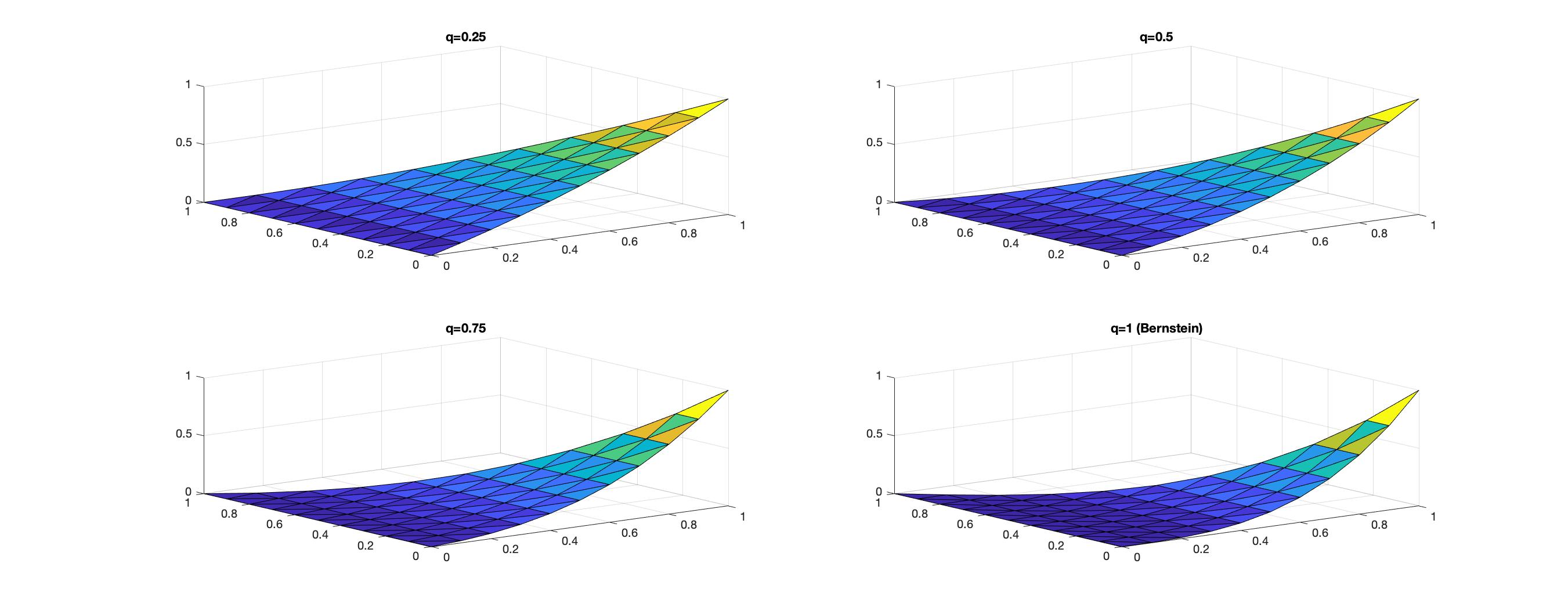} 
\caption{The effect of changing the q-parameter on the cubic $q$--Bernstein polynomial $B_{003}^3$}\label{fig:b003}
\end{center}
\end{figure}

\begin{figure}[!h]
\begin{center}
\includegraphics[width=\textwidth]{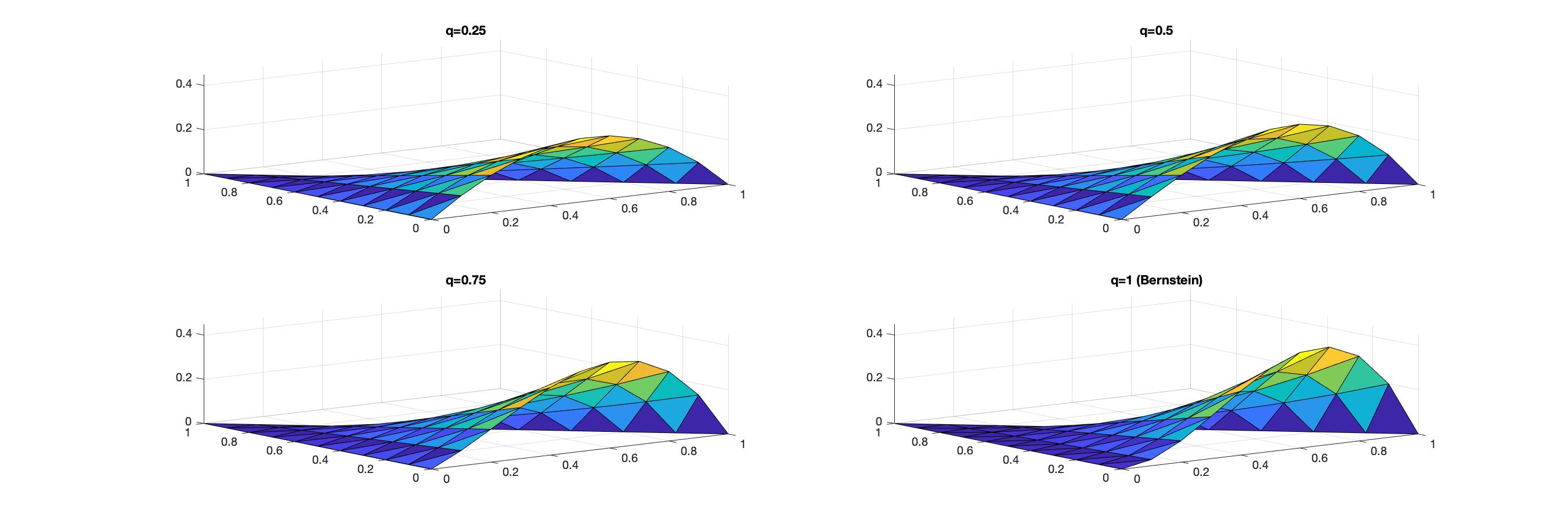} 
\caption{The effect of changing the q-parameter on the cubic $q$--Bernstein polynomial $B_{012}^3$}\label{fig:b012}
\end{center}
\end{figure}

\begin{figure}[!h]
\begin{center}
\includegraphics[width=\textwidth]{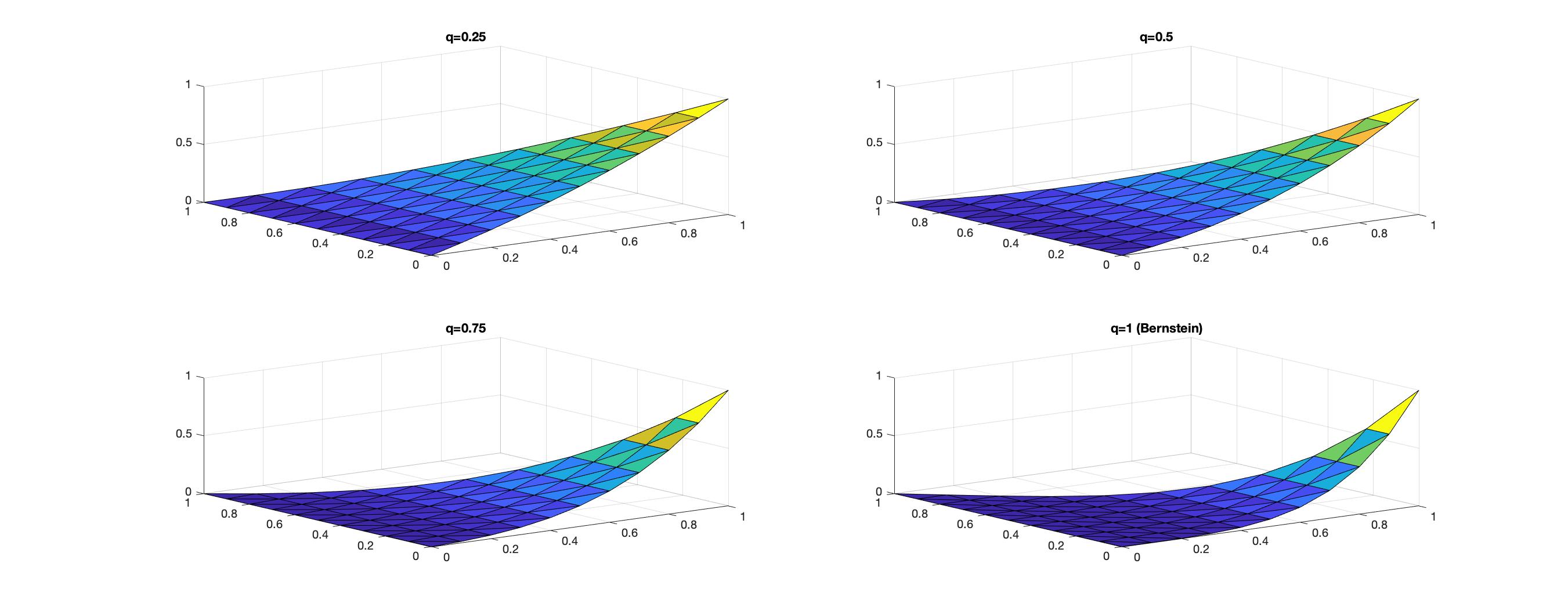} 
\caption{The effect of changing the q-parameter on the cubic $q$--Bernstein polynomial $B_{004}^4$}\label{fig:b004}
\end{center}
\end{figure}

Figures \ref{fig:b003} and \ref{fig:b012} show the graphs of some cubic $q$--Bernstein basis functions for different values of the shape parameter $q$. Let us recall that for $q=1$ the $q$-analogues of the triangular Bernstein
polynomials coincide with the usual triangular Bernstein polynomials. 
Figure \ref{fig:b004} shows the graph of the quartic $q$--Bernstein basis  function $B_{004}^4$ for different values of $q$. So, in those figures
we can see the effect of changing the parameter $q$.

\section{An evaluation algorithm of the $q$--B\'ezier representation of polynomials over triangles}
\label{sec:eval.Alg}

Let us consider the polynomials $\gamma^n(u,v)$  that are linear combinations of  the polynomials in \eqref{def.qBer}. Then we say that $\gamma^n(u,v)$ is a $q$--B\'ezier representation of a polynomial of degree $n$ and it is given by
\begin{equation}\label{def.bpol}
\gamma^n(u,v):=\sum_{i+j+k=n}b_{ijk}^nB_{ijk}^n(u,v),
\end{equation}
where $b_{ijk}^n\in \mathbb{R}$.
We are going to see that polynomials given by the representation \eqref{def.bpol} can be evaluated using a de Casteljau type algorithm. Our proposed method is given by Algorithm \ref{alg:eval}. This algorithm consists on $n$ steps and it produces the evaluation of $\gamma^n(u,v)$ at any point $P=uT_1+vT_2+wT_3$. If the point is in $\mathcal{T}$, then every step of the algorithm is formed by linear convex combinations of the previously computed quantities. In fact, for $q=1$ it coincides with the classical de Casteljau algorithm for the triangular Bernstein polynomials.
\begin{algorithm}[!h]
 \caption{Evaluation of the polynomial $\gamma^n(u,v)$}\label{alg:eval}
 \begin{algorithmic}
  \Require $b_{ijk}^n$ for all $i+j+k=n$, $u$, $v$ and $q\in (0,1]$.
  \Ensure $f_{000}^{(n)}$ (the evaluation of $\gamma^n$ on $P=uT_1 +vT_2+ wT_3$)
  \State $f_{ijk}^{(0)}=b_{ijk}^n$ 
  \For{$r = 1:n$}
  \ForAll{$i+j+k =n-r $}
        	        	\State  $f_{ijk}^{(r)}=q^{k}uf_{i+1,j,k}^{(r-1)}+q^{k}vf_{i,j+1,k}^{(r-1)}+(1-q^{k}u-q^{k}v)f_{i,j,k+1}^{(r-1)}$
    	\EndFor
	\EndFor
 \end{algorithmic}
\end{algorithm}

\begin{remark}
 An evaluation algorithm formed by linear convex combinations for the $q$-Bernstein polynomials of one variable can be obtained as a particular case of Algorithm \ref{alg:eval}. In fact, if we only consider points such that $v=0$ (or analogously, $u=0$), Algorithm \ref{alg:eval} produces the evaluation of a q-Bernstein polynomial on the line defined by $T_1$ and $T_3$ (respectively, $T_2$ and $T_3$).
\end{remark}

In order to show that Algorithm \ref{alg:eval} gives the evaluation of the polynomial given by \eqref{def.bpol} at the point $P$, in Proposition \ref{prop.eva} we give an explicit formula for the intermediate quantities $f_{ijk}^{(r)}$ computed at every step of the algorithm. 

\begin{proposition}\label{prop.eva}
For $0\leq r \leq n$ and for all $i+j+k=n-r$, let us define the quantities by a recursive relationship:
\begin{align}
f_{ijk}^{(0)}&=b_{ijk}^n,\nonumber\\
f_{ijk}^{(r)}&=q^kuf_{i+1,j,k}^{(r-1)}+q^kvf_{i,j+1,k}^{(r-1)}+(1-q^ku-q^kv)f_{i,j,k+1}^{(r-1)}.\label{f.eval}
\end{align}
Then we have that $\gamma^n(u,v)=f_{000}^n$.
\end{proposition}
\begin{proof}
Let us start by using \eqref{recq} to deduce a formula for the evaluation of the polynomial $\gamma^n(u,v)$ in terms of the $q$-Bernstein polynomials of degree $n-1$,
\begin{align}\label{aux}
\gamma^n(u,v)&=\sum_{i+j+k=n}f_{ijk}^{(0)}B_{ijk}^n(u,v)\nonumber\\
&=\sum_{i+j+k=n}f_{ijk}^{(0)}(q^k u B_{i-1,j,k}^{n-1}(u,v) + q^k vB_{i,j-1,k}^{n-1}(u,v) + (1-q^{k-1}u-q^{k-1}v)B_{i,j,k-1}^{n-1}(u,v)).
\end{align}
By rearranging the summands in \eqref{aux}, we see that
\begin{align*}
\gamma^n(u,v)&=\sum_{i+j+k=n-1}(q^ku f_{i+1,j,k}^{(0)} + q^kvf_{i,j+1,k}^{(0)}+ (1-q^ku-q^kv)  f_{i,j,k+1}^{(0)}) B_{ijk}^{n-1}(u,v)\\
&= \sum_{i+j+k=n-1} f_{ijk}^{(1)}B_{ijk}^{n-1}(u,v),
\end{align*}
where $f_{ijk}^{(1)}$ is given by \eqref{f.eval}. In general, we can apply this strategy $r$ times to deduce the following  expression for the evaluation of $\gamma^n(u,v)$ in terms of the polynomials $B^{n-r}_{ijk}(u,v)$:
\begin{align*}
\gamma^n(u,v)&= \sum_{i+j+k=n-r} f_{ijk}^{(r)}B_{ijk}^{n-r}(u,v).
\end{align*}
In particular,  for $r=n$ we have that
$\gamma^n(u,v)=f_{000}^{(n)}B_{000}(u,v)=f_{000}^{(n)}$. Hence, Algorithm \ref{alg:eval} gives the evaluation of the polynomial $\gamma^n$ of degree $n$ given by  \eqref{def.bpol} at $(u,v)$.
\end{proof}

 We have that the triangular $q$--Bernstein polynomials form a partition of unity as a direct consequence of the evaluation algorithm.

\begin{corollary}\label{prop.unity}
The triangular $q$--Bernstein polynomials defined by \eqref{def.qBer} form a  partition of unity, i.e.,

\begin{equation}\label{part.unity}
\sum_{i+j+k=n} B_{ijk}^n (u,v)=1.
\end{equation}
\end{corollary}
\begin{proof}
Let us define a  polynomial  $\gamma^n(u,v)$ following \eqref{def.bpol} whose coefficients $b_{ijk}^n$ are all $1$'s. Since Algorithm \ref{alg:eval} evaluates  $\gamma^n(u,v)$ at every point inside $\mathcal{T}$ through consecutive combinations of coefficients summing up to 1  that take $b_{ijk}^n$ as its initial input, we have that $f_{ijk}^{(r)}=1$ for all $i,j,k,r$. Hence, $\gamma^n(u,v)=f_{000}^{(n)}=1$ for every $u,v\geq 0$ such that $u+v\leq 1$ and we can conclude that $\sum_{i+j+k=n} B_{ijk}^n(u,v)= \gamma^n(u,v)=1$.
\end{proof}

\section{$q$--Bernstein basis in a triangle and $q$--B\'ezier patches}

The following result shows that it is possible to write the polynomial $\gamma^n(u,v)$ given by \eqref{def.bpol} as a linear combination of the $q$--Bernstein basis of degree $n+1$, giving a degree elevation formula.

\begin{proposition} \label{degree.raising}
The polynomial $\gamma^n(u,v)$ of degree $n$ given by \eqref{def.bpol} can be written as a linear combination of the triangular $q$--Bernstein polynomials of degree $n+1$,

\[\gamma^{n}(u,v)=\sum_{i+j+k=n+1}b_{ijk}^{n+1}B_{ijk}^{n+1}(u,v),\]
where 
\[b_{ijk}^{n+1}=\frac{[n+1-k]}{[n+1]}\frac{i}{i+j}q^kb_{i-1,j,k}^{n}+\frac{[n+1-k]}{[n+1]}\frac{j}{i+j}q^kb_{i,j-1,k}^{n}+\frac{[k]}{[n+1]}b_{i,j,k-1}^{n}.\]
\end{proposition}
\begin{proof}
Let us start by writing $\gamma^n(u,v)$ as a linear combination of the triangular $q$--Bernstein polynomials of degree $n$ following \eqref{def.bpol},
\begin{align*}
\gamma^n(u,v)=& \sum_{i+j+k=n}b_{ijk}^nB_{ijk}^n(u,v)=\sum_{i+j+k=n}b_{ijk}^nB_{ijk}^n(u,v)(q^{k}u +q^{k}v+ 1 -q^{k}u-q^{k}v).
\end{align*}
Using the definition of $B_{ijk}(u,v)$ given by \eqref{def.qBer}, we see that 
\begin{align}
\gamma^n(u,v)=& \sum_{i+j+k=n}q^k b_{ijk}^n\qbinom{n}{k}\binom{i+j}{i}u^{i+1}v^j\prod_{s=0}^{k-1}(1-q^su-q^sv)\nonumber\\
&+ \sum_{i+j+k=n}q^k b_{ijk}^n\qbinom{n}{k}\binom{i+j}{i}u^{i}v^{j+1}\prod_{s=0}^{k-1}(1-q^su-q^sv)\nonumber \\
&+ \sum_{i+j+k=n} b_{ijk}^n\qbinom{n}{k}\binom{i+j}{i}u^{i}v^{j}\prod_{s=0}^{k}(1-q^su-q^sv).\label{aa}
\end{align}

Let us now choose indices $i$, $j$ and $k$ such that $i+j+k=n+1$ and let us see what are the coefficients for $B_{ijk}^{n+1}(u,v)$ in \eqref{aa},
\begin{align*}
q^k b_{i-1,j,k}^n\qbinom{n}{k}\binom{i+j-1}{i-1}u^{i}v^j\prod_{s=0}^{k-1}(1-q^su-q^sv)\\
+ q^k b_{i,j-1,k}^n\qbinom{n}{k}\binom{i+j-1}{i}u^{i}v^j\prod_{s=0}^{k-1}(1-q^su-q^sv)\\
 + b_{i,j,k-1}^n\qbinom{n}{k-1}\binom{i+j}{i}u^{i}v^{j}\prod_{s=0}^{k-1}(1-q^su-q^sv).
\end{align*}
Using the definition of $B_{ijk}^{n+1}(u,v)$, \eqref{def.qBer}, we can rewrite the last expression by \eqref{def.qbi} as

\begin{align*}
\left(\frac{[n+1-k]}{[n+1]}\frac{i}{i+j}q^kb_{i-1,j,k}^{n}+\frac{[n+1-k]}{[n+1]}\frac{j}{i+j}q^kb_{i,j-1,k}^{n}+\frac{[k]}{[n+1]}b_{i,j,k-1}^{n}\right)B_{ijk}^{n+1}&=b_{ijk}^{n+1}B_{ijk}^{n+1}.      
\end{align*}
Hence, 
\[\gamma^{n}(u,v)=\sum_{i+j+k=n+1}\left(\frac{[n+1-k]}{[n+1]}\frac{i}{i+j}q^kb_{i-1,j,k}^{n}+\frac{[n+1-k]}{[n+1]}\frac{j}{i+j}q^kb_{i,j-1,k}^{n}+\frac{[k]}{[n+1]}b_{i,j,k-1}^{n}\right)B_{ijk}^{n+1}(u,v),\]
and the statement of the proposition holds.
\end{proof}
Using the degree elevation property given by Proposition \ref{degree.raising}, we now prove that the triangular $q$--Bernstein polynomials of degree $n$ form a basis of the space $\mathbb{P}^n$ of bivariate polynomials of degree less than or equal to $n$.
\begin{theorem}
The set
\begin{equation}\label{basis}
\mathcal{B}^n_q=\{B_{ijk}^n(u,v)\}_{i+j+k=n}
\end{equation}
of triangular $q$--Bernstein polynomials of degree $n$ is a basis for the space $\mathbb{P}^n$.
\end{theorem}

\begin{proof}
The dimension of the space of bivariate polynomials $\mathbb{P}^n$ is equal to the number of triangular $q$--Bernstein polynomials of degree $n$, $\binom{n+2}{2}$ (see page 1 of \cite{LS}). Hence, we can prove that the triangular $q$--Bernstein polynomials form a basis of $\mathbb{P}^n$ by showing that  the monomial basis $\{x^\alpha y^\beta \}_{0\leq \alpha + \beta \leq n}$  is in the span of $\mathcal{B}^n_q$. By Corollary \ref{prop.unity}, we know that $1\in\langle\mathcal{B}^n_q\rangle$.

Let us proceed by induction over $n$.  For $n=1$, we have that $1=u+v+(1-u-v)$ and that, by the definition of the barycentric coordinates \eqref{coord}, $x=x_1 u + x_2 v + x_3 (1-u-v)$ and $y=y_1 u + y_2 v + y_3 (1-u-v)$.
Now let us assume that the theorem holds for the polynomials of degree $n$. Then, $x^\alpha y^\beta \in\langle\mathcal{B}^n_q\rangle$ for all $0\leq \alpha + \beta \leq n$ and, by the degree elevation property given by Proposition \ref{degree.raising}, $x^\alpha y^\beta \in\langle\mathcal{B}^{n+1}_q\rangle$ for all $0\leq \alpha + \beta \leq n$. 
Hence, let us choose $\alpha$ and $\beta$  such that $\alpha +\beta =n+1$. Without loss of generality, let us assume that $\alpha \geq 1$. By the induction hypothesis  we know that $x^{\alpha-1}y^\beta=\sum_{i+j+k=n}b_{ijk}^nB_{ijk}^n(u,v)$. Hence,
\begin{align*}
x^{\alpha}y^\beta&=(ux_1+vx_2+(1-u-v)x_3)x^{\alpha-1}y^\beta\\
&=\sum_{i+j+k=n}(ux_1+vx_2+(1-u-v)x_3) b_{ijk}^nB_{ijk}^n(u,v)\\
&=\sum_{i+j+k=n}\left((x_1-x_3+q^k x_3)u+(x_2-x_3+q^k x_3)v+x_3(1-q^ku-q^kv) \right) b_{ijk}^nB_{ijk}^n(u,v).
\end{align*}
In the previous formula, we write  $B^n_{ijk}(u,v)$ in terms of the triangular $q$--Bernstein polynomials of degree $n+1$ (see Proposition \ref{degree.raising}),
\begin{align*}
x^{\alpha}y^\beta&=\sum_{i+j+k=n} \frac{[n+1-k]}{[n+1]}\frac{i+1}{i+j+1}(x_1-x_3+q^k x_3)b_{ijk}^{n}B_{i+1,j,k}^{n+1}(u,v)\\
&+\sum_{i+j+k=n} \frac{[n+1-k]}{[n+1]}\frac{j+1}{i+j+1}(x_2-x_3+q^k x_3)b_{ijk}^{n}B_{i,j+1,k}^{n+1}(u,v)\\
&+\sum_{i+j+k=n}\frac{[k+1]}{[n+1]}x_3b_{ijk}^{n}B_{i,j,k+1}^{n+1}(u,v).
\end{align*}
Finally, by rearranging the terms of the sum, we can collect the coefficients of each triangular $q$--Bernstein polynomial and rewrite $x^{\alpha}y^\beta$ as
\begin{align*}
x^{\alpha}y^\beta=\sum_{i+j+k=n+1}& \left(\frac{[n+1-k]}{[n+1]}\frac{i}{i+j}(x_1-x_3+q^k x_3)b_{i-1,jk}^{n} + \frac{[n+1-k]}{[n+1]}\frac{j}{i+j}(x_2-x_3+q^k x_3)b_{i,j-1,k}^{n}\right.\\
&\left.+\frac{[k]}{[n+1]}x_3b_{i,j,k-1}^{n}\right)B_{ijk}^{n+1}(u,v).
\end{align*}
And so, $x^{\alpha}y^\beta \in \langle\mathcal{B}^{n+1}_q\rangle$ and the result holds for $n+1$.
\end{proof}

In \cite{TPQB}, $q$--B\'ezier patches on a rectangle were introduced by using tensor products of univariate $q$--Bernstein polynomials. Now we extend the notion of  parametric $q$--B\'ezier curves to parametric $q$--B\'ezier triangular surfaces. We define the {\it $q$--B\'ezier patch} of degree $n$ as the parametric surface given by 
\begin{equation}\label{def.bpatch}
Q(u,v):=\sum_{i+j+k=n}P_{ijk}^nB_{ijk}^n(u,v),
\end{equation}
where $(u,v,w=1-u-v)$ are the barycentric coordinates of a point $P=(x,y)$ with respect to the triangle $T$ (see (\ref{coord})) and 
$P_{ijk}^n\in \mathbb{R}^3$ are the {\it control points} of the $q$--B\'ezier patch.

Observe that we can apply componentwise the evaluation algorithm given by Algorithm 1 to derive a de Casteljau type evaluation algorithm for the $q$--B\'ezier patch.

\begin{figure}[!h]
\begin{center}
\includegraphics[width=\textwidth]{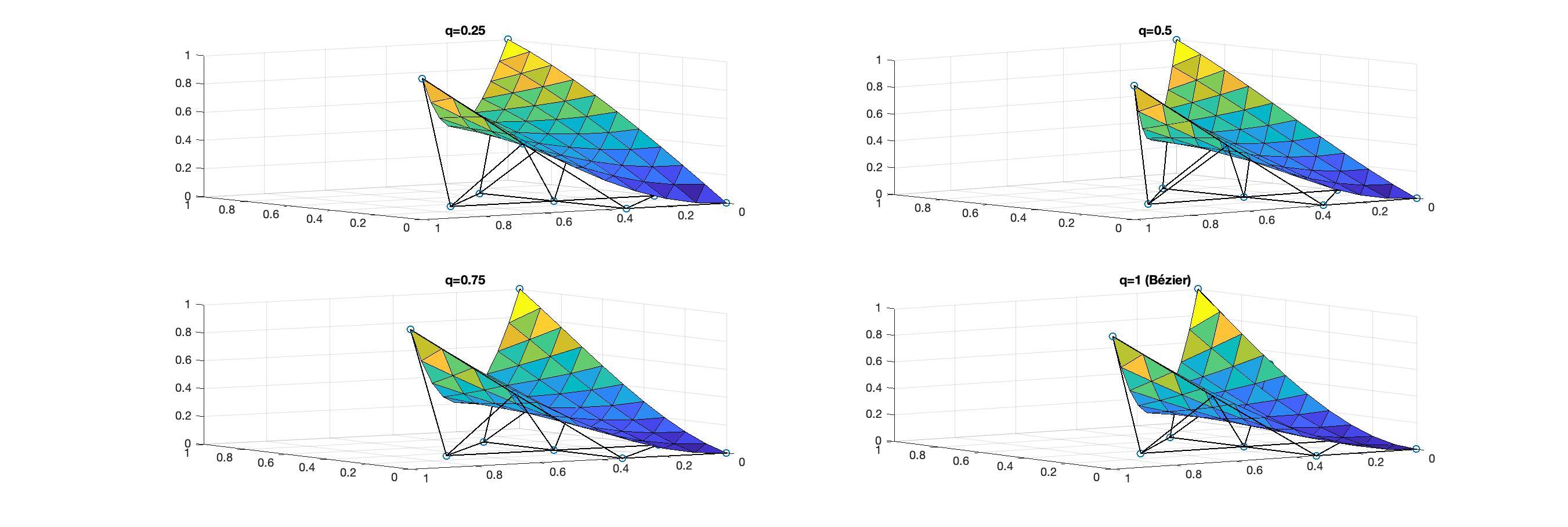} 
\caption{Cubic q-B\'ezier patch}\label{fig:p2}
\end{center}
\end{figure}

\begin{figure}[!h]
\begin{center}
\includegraphics[width=\textwidth]{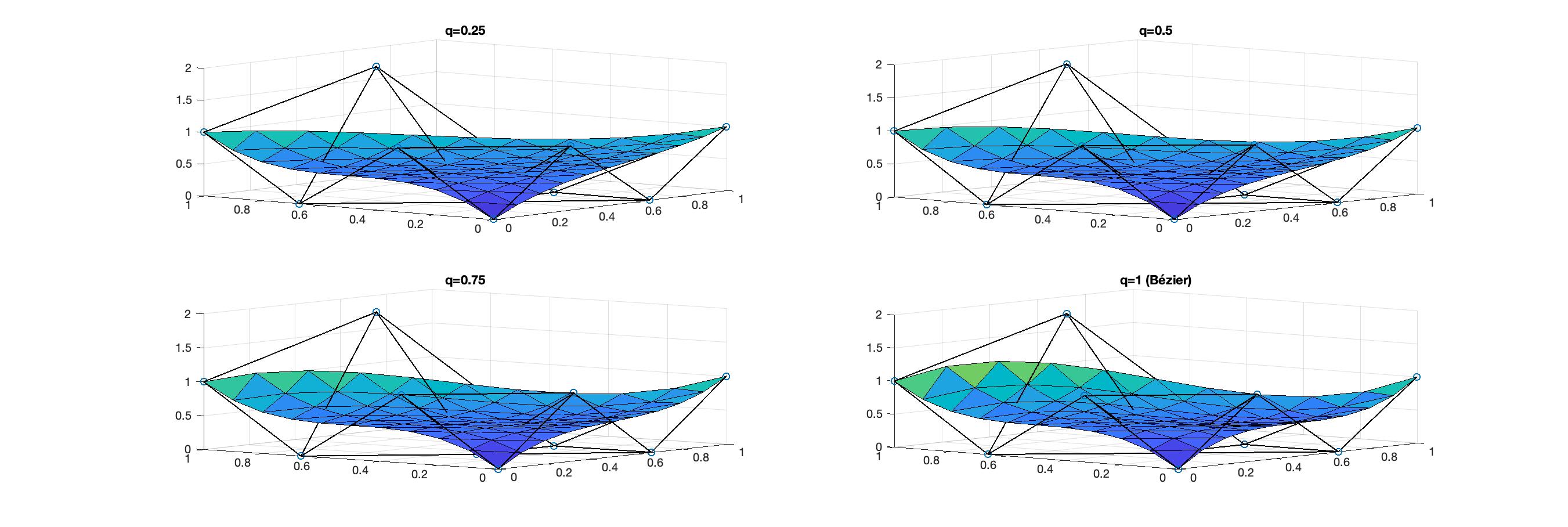} 
\caption{Cubic q-B\'ezier patch}\label{fig:p3}
\end{center}
\end{figure}

Figure \ref{fig:p2} shows cubic $q$--B\'ezier patch with control 
points $P_{300}^3=(0,0,0)$, $P_{210}^3=(0,1/3,0)$,
$P_{120}^3=(0,2/3,1/2)$, $P_{030}^3=(0,1,1)$,
$P_{201}^3=(1/3,0,0)$, $P_{111}^3=(1/3,1/3,0)$,
$P_{021}^3=(1/3,2/3,0)$,
$P_{102}^3=(2/3,0,1/2)$, $P_{012}^3=(2/3,1/3,0)$ and
$P_{003}^3=(1,0,1)$
over a triangle  for different values of the shape parameter $q$. Figure \ref{fig:p3} shows another cubic $q$--B\'ezier patch with control points
$P_{300}^3=(0,0,0)$, $P_{210}^3=(0,1/3,1)$,
$P_{120}^3=(0,2/3,0)$, $P_{030}^3=(0,1,1)$,
$P_{201}^3=(1/3,0,1)$, $P_{111}^3=(1/3,1/3,0)$,
$P_{021}^3=(1/3,2/3,2)$
$P_{102}^3=(2/3,0,0)$, $P_{012}^3=(2/3,1/3,0)$ and
$P_{003}^3=(1,0,1)$
 over a triangle for different values of $q$.

\section{Stabiliy of the basis}

In Section \ref{sec:eval.Alg}, an algorithm for the evaluation of polynomials over triangles using the $q$-Bézier representation has
been presented. Since all its steps are linear convex combinations,
the evaluation algorithm is stable. The accuracy of the evaluations provided by the algorithm depends on
the stability of the algorithm and on the conditioning of the evaluated polynomial. In turn, this conditioning 
depends on the basis used to represent the function. 
In \cite{L-P}, it was proved that the triangular Bernstein bases are  best conditioned among all nonnegative 
bases of the corresponding spaces of multivariate polynomials in the sense that there is not a better conditioned nonnegative basis. In \cite{DP.RTS} it was also proved the same fact
for the rational triangular Bernstein bases, a generalization of those triangular Bernstein bases.
Now the conditioning of polynomials represented in the q-Bernstein basis and the usual Bernstein basis over a triangle is going to be compared.

First, let us recall the definition of relative condition number (see \cite{L-P}). For any $x\in \Omega\subset \mathbb{R}^s$, basis
$u=(u_0,\ldots,u_n)$ and $f=\sum_{i=0}^n c_i u_i$,
let us consider the relative condition number for the 
evaluation of $f$ at $x$:
\[
cond(u;f,x)=\frac{\sum_{i=0}^n\left| c_i u_i(x)\right|}
	{\left\|f \right\|_{\infty}}=\frac{\sum_{i=0}^n\left| c_i u_i(x)\right|}
	{\left\|\sum_{i=0}^n c_i u_i\right\|_{\infty}}
\]

The following result will be an important tool when comparing the conditioning of a polynomial for the two mentioned
representations.
\begin{lemma} \label{lem:cond}(cf. Lemma 3.1 of \cite{L-P})
Let  $\mathcal{U}$ be a finite dimensional vector space of functions
defined on $\Omega\subset \mathbb{R}^s$. Let $u,v$ be two bases
of nonnegative functions of $\mathcal{U}$. Then
\[
	cond(u;f,x)\leq cond(v;f,x),\quad \forall\,f\in\mathcal{U},\ \forall x\in\Omega,
\]
if and only if the matrix $A$ such that $v=uA$ is nonnegative.
\end{lemma}
In order to determine that all the entries  of the matrix of change of bases from the q-Bernstein basis to the 
Bernstein basis are nonnegative the following 
auxiliary result will be necessary.
\begin{lemma}\label{lem:aux}
For $0<q<1$ and $r\ge 1$, it is satisfied that
\begin{equation}\label{eq:res.est}
	\prod_{s=0}^{r-1} (1-q^su-q^sv)=\sum_{i+j+k=r} c_{ijk}^r u^iv^j(1-u-v)^k
\end{equation}
with $c_{ijk}^r\ge 0$ for all $i,j,k\ge 0$ with $i+j+k=r$.
\end{lemma}
\begin{proof}
Let us prove the result by induction on $r\ge 1$. For $r=1$, Equation \eqref{eq:res.est}
is satisfied
with $c_{100}^1=c_{010}^1=0$ and $c_{001}^1=1$, and the result holds. Let us assume
that the result satisfies for $r=t\ge 1$ and let us prove it for $r=t+1$. So, for $r=t+1$ 
formula \eqref{eq:res.est}  can be written as
\[
	\prod_{s=0}^{t} (1-q^su-q^sv)=
	\left[\prod_{s=0}^{t-1} (1-q^su-q^sv)\right](1-q^tu-q^tv).
\]
By the induction hypothesis and operating we have that
\begin{align*}
	\prod_{s=0}^{t} (1-q^t u-q^tv)&=
	\left(\sum_{i+j+k=t} c_{ijk}^t u^iv^j(1-u-v)^k\right) (1-q^tu-q^tv) \\
	&=\left(\sum_{i+j+k=t} c_{ijk}^t u^iv^j(1-u-v)^k\right) \left[(1-u-v)+(1-q^t)u+(1-q^t)v\right]
\end{align*}
with $c_{ijk}^t\ge 0$ for all $i+j+k=t$. Then, we can deduce that
\begin{align}\label{eq:interm}
	\prod_{s=0}^{t} (1-q^su-q^sv)
	&=\left(\sum_{i+j+k=t} c_{ijk}^t u^iv^j(1-u-v)^{k+1}\right) +
	\left(\sum_{i+j+k=t} (1-q^t)c_{ijk}^t u^{i+1}v^j(1-u-v)^k\right)  \nonumber \\
	&\ \ +\left(\sum_{i+j+k=t} (1-q^t)c_{ijk}^t u^iv^{j+1}(1-u-v)^k\right).
\end{align}
Taking into account that for $0<q<1$, $1-q^t>0$ holds. Then, since $1-q^t>0$ and $c_{ijk}^t\ge 0$ 
formula \eqref{eq:interm} can written as
\[
\prod_{s=0}^{t} (1-q^su-q^sv)=\sum_{i+j+k=t} c_{ijk}^{t+1} u^iv^j(1-u-v)^k
\]
with $c_{ijk}^{t+1}\ge 0$ and the induction holds.
\end{proof}
In the following result we prove that the matrix of change between the  $q$-Bernstein Bernstein basis
and the usual Bernstein basis over a triangle has all its entries nonnegative and, as a consequence,
the usual Bernstein basis is better conditioned than the $q$-Bernstein basis.
\begin{theorem}
Let $b_1=(b_{ijk}^n(u,v))_{i+j+k=n}$ be the Bernstein triangular polynomials of degree $n$ given by 
\[
	b_{ijk}^n(u,v)=\frac{n!}{i!\cdot j!\cdot k!}u^iv^j(1-u-v)^{n-i-j},\quad i,j,k\ge 0\text{ with } i+j+k=n,
\]
and let $b_2=(B_{ijk}^n(u,v))_{i+j+k=n}$ be the q-analogues of the triangular Bernstein polynomials for $0<q<1$. Then,
\begin{itemize}
\item[i.] all the entries of the matrix $A$ of change of bases such that $b_2=b_1A$ are nonnegative, and
\item[ii.] $cond(b_1;f,(u,v))\leq cond(b_2;f,(u,v)),\quad \forall\,f\in\left\langle b_1\right\rangle,\forall u,v\ge 0\text{ such that }u+v\leq 1$.
\end{itemize}
\end{theorem}
\begin{proof}
\begin{itemize}
	\item[i.] By applying Lemma \ref{lem:aux} in formula \eqref{def.qBer} it is deduced that
	\begin{align*}
		B_{ijk}^n(u,v)&=\qbinom{n}{k}\binom{i+j}{i}u^iv^j\prod_{s=0}^{k-1}(1-q^su-q^sv) \\
			&=\qbinom{n}{k}\binom{i+j}{i}u^iv^j\sum_{r+s+t=k} c_{rst}^k u^rv^s(1-u-v)^t \\
			&=\sum_{r+s+t=k} c_{rst}^k\qbinom{n}{k}\binom{i+j}{i}u^{i+r}v^{j+s}(1-u-v)^{t},\quad i+j+k=n,
	\end{align*}
with $c_{rst}^k\ge 0$ for all $r,s,t\ge 0$ with $r+s+t=k$. Taking into account that $(i+r)+(j+s)+t=n$, that $c_{rst}^k\ge 0$,
and that $\qbinom{n}{k},\binom{i+j}{i}>0$, we can conclude that
\[
B_{ijk}^n(u,v)=\sum_{i+j+k=n}a_{ijk}^n b_{ijk}^n(u,v)
\]
where $a_{ijk}^n\ge 0$ for all $i+j+k=n$ . Hence i follows.
\item[ii.] It is an straightforward consequence of part i and Lemma \ref{lem:cond}.
\end{itemize} 
\end{proof}

\section{Conclusions and future work}

We have defined and analyzed $q$-Bernstein basis functions over a triangular domain, which are a natural extension of the univariate ones and an alternative to the extension of \cite{LWAG}. In contrast to \cite{LWAG}, our evaluation algorithm is always formed by linear convex combinations for every point in the triangle, which implies nice stability and geometric properties. We have seen that the matrix of change of basis from the Bernstein basis to the $q$-Bernstein basis is nonnegative, which implies that the Bernstein basis is better conditioned. Our new basis functions also satisfy some recurrence relations, form a partition of unity and are a basis for the space of polynomials of total degree less than or equal to $n$ on a triangle. In addition to a de Casteljau type evaluation algorithm, our basis functions also satisfy other convenient properties for design purposes such as degree elevation.

As a future work, we shall consider the possible extension of the subdivision property to the new representation of a triangular polynomial. An adequate extension of the construction of the presented basis to a general multivariate case will also be analyzed.

\section*{Acknowledgments}

We would like to thank the anonymous referees for their valuable insight and suggestions.

\end{document}